\documentclass[reqno]{amsart}
\theoremstyle{plain}
\usepackage{graphicx}
\newtheorem{thm}{\it Theorem}[section]
\newtheorem{prop}[thm]{\it Proposition}
\newtheorem{cor}[thm]{\it Corollary}
\newtheorem{lem}[thm]{\it Lemma}
\theoremstyle{remark}
\newtheorem{defn}[thm]{Def{}inition}
\newtheorem{rem}[thm]{Remark}
\newtheorem{exa}[thm]{Example}
\numberwithin{equation}{section}
\usepackage{enumerate}

\allowdisplaybreaks
\begin{document}
	\baselineskip08pt
	\title [Frames generated by graphs]{Frames generated by graphs}

	\author[Deepshikha]{ Deepshikha }
	\address{Deepshikha, Department of Mathematics, Shyampur Siddheswari Mahavidyalaya, University of  Calcutta-711312,  India.}
	\email{dpmmehra@gmail.com}

	\begin{abstract}
		Frames are the most natural generalization of orthonormal bases that allow the inclusion of redundant systems. In this article, we introduce the concept of frames generated by graphs in finite-dimensional spaces and study their properties. Let $G$ be a simple graph of $n$ vertices with Laplacian matrix $L$. We define the notions of $G(n,k)$-frames and $L_G(n,k)$-frames associated with the graph $G$. We obtain the family of dual frames of $L_G(n,k)$-frames and $G(n,k)$-frames. It is shown that non-regular graphs cannot generate tight frames. Then we establish a characterization of tight $G(n,k)$-frames in terms of the adjacency spectra of regular graphs. Besides, we provide a frame theoretic proof of an existing graph property. Finally, we show that one can use complete graphs to generate tight frames.
	\end{abstract}
	
	\subjclass[2020]{42C15;  42C40; 46C05; 05C50}
	
	\keywords{  Finite frames; $G(n,k)$-frames; $L_G(n,k)$-frames; Tight frames; Regular graphs; Laplacian matrix.}
	
	\maketitle
	
	\baselineskip15pt
	
	\section{Introduction}

	Frame, graph, and matrix theories are the three applicable, useful, and fundamental branches of mathematics. They have their own deep significance in their domain of studies. However, it is always interesting to discover some relation among different fields of mathematics. In this article, we establish a connection among these three branches of mathematics and provide several results. Before discussing the work presented in this paper, we give a brief introduction to spectral graph theory (i.e., the theory of matrices in connection with graphs) and frame theory in the next two paragraphs. \vspace{4pt}

		Graph theory is one of the oldest branches of mathematics. Graphs are used to model many real-life problems. Spectral graph theory studies the properties of graphs in terms of the associated matrix properties and vice-versa. Various matrices can be associated with a graph, like adjacency matrix, Laplacian matrix, distance matrix, etc. Along with the deep theoretical aspect of spectral graph theory, it has many real-life applications, which include measuring the connectivity of communication networks, hierarchical clustering, ranking of hyperlinks in search engines, mixing time of a random walk on graphs, image segmentation, and so on, see \cite{AO,B,BH,CS}. 

	\vspace{4pt}
	
	The concept of frames was given by Duffin and Schaeffer \cite{DS} and revived by Daubechies, Grossmann, and Meyer \cite{DGM}, who showed its importance in data processing. Frame is the redundant sequence of vectors that span a vector space and provides infinitely many representations of a vector.  This redundancy is key to the resilience of frames to various disturbances like noises, erasures, quantization, etc., see \cite{BP, BP1,  CK, HP, SH}.  Apart from the application of frames in numerous practical problems like signal and image processing, sparse representations, data compression, etc., frames are also useful in many theoretical problems of analysis. For instance, a very famous Kadison-Singer problem has some equivalent formulations in frame theory, and one of them is known as Weaver’s conjecture. Recently, Marcus et al. \cite{MSS} solved the Kadison-Singer problem by providing a positive answer to the Weaver’s conjecture.  \vspace{4pt}

	Signal transmission and processing using frames involves the decomposition of a signal in terms of frame coefficients and its reconstruction using dual frames, see \cite{CK1, C2, LH1, PHM}. In this paper, we generate finite frames using the Laplacian matrix of a graph. Then we obtain a result, which provides the family of all dual frames of the frames generated by graphs. \vspace{4pt}

	Erasure is a common problem in wireless sensor networks due to transmission losses or the occurrence of noises during transmission, see \cite{NAS, BP, DA, HP}. Uniform tight frames are known to be optimal for handling up to one erasure, see \cite{BP, CK}. We characterize tight frames generated by graphs in terms of graph-theoretic properties.   We show that one can construct tight frames using regular graphs with two eigenvalues, and tight frames generated by graphs are always uniform frames. Hence, we can use graphs to construct uniform tight frames. Many authors studied frame theoretic properties using matrices associated with graphs. For instance, the existence of two-uniform frames (or Grassmannian frames), which are optimal for one erasure and two erasures, are studied by using the Seidel adjacency matrix of graphs with two eigenvalues, see \cite{BP, SH}. In this article, we show the frame theoretic proof of an existing graph property.
	
	This paper is organized as follows. Section 2 provides the basics of finite frame theory and graph theory required in the rest of the article. In section 3, we show how to generate a frame using the Laplacian matrix $L$ of a graph $G$. We introduce the concept of $L_G(n,k)$-frame and $G(n,k)$-frame associated with the graph $G$. Then we give an example of a $G(n,k)$-frame, which is not an $L_G(n,k)$-frame. We show that the frames generated by a graph are unique up to unitary equivalence. We also obtain the family of canonical dual and alternate dual frames of any $G(n,k)$-frame and $L_G(n,k)$-frame. In Section 4, we study a connection between tight frames and regular graphs. We provide a necessary condition of tight $G(n,k)$-frames in terms of regular graphs. Further, we characterize tight frames generated by graphs. At the end of the article, we provide a frame theoretic proof of one of the known graph property.

	\section{Preliminaries}
	\subsection{Frames}
	Throughout the paper, we consider finite frames on $k$-dimensional Hilbert space $\mathbb{C}^k$ and $[n]=\{1,2,\ldots,n\}$. We denote a zero matrix or zero vector of appropriate order by $\mathbf{0}$. We use $\{e_i\}_{i\in[n]}$ to denote the standard canonical orthonormal basis of the Hilbert space $\mathbb{C}^n$. If $\{ f_i\}_{i\in[n]}$ is a sequence of vectors  in $\mathbb{C}^k$, then it is called a \emph{frame}  (or \emph{finite frame}) for  $\mathbb{C}^k$ if there exist positive constants $A$ and $B$ such that for any $f\in\mathbb{C}^k$, we have
	\begin{align}\label{eq2.1}
		A\|f\|^2\leq\sum_{i\in[n]}|\langle f, f_i\rangle|^2\leq B\|f\|^2.
	\end{align}
	The inequality \eqref{eq2.1} is called the \emph{frame inequality}. The constants $A$ and $B$ are called the \emph{lower frame bound} and \emph{upper frame bound} of $\{f_i\}_{i\in[n]}$, respectively. If the upper inequality holds in \eqref{eq2.1}, then $\{f_i\}_{i\in[n]}$ is called a \emph{Bessel sequence}, and $B$ is called a \emph{Bessel bound}. If the constants $A=B$, then $\{f_i\}_{i\in[n]}$ is called a \emph{tight frame}, and in this case, $\{f_i\}_{i\in[n]}$ is also said to be an \emph{$A$-tight frame}. For $A=B=1$, $\{f_i\}_{i\in[n]}$ is called a \emph{Parseval frame}. If $\{f_i\}_{i\in[n]}$ is a frame and there exists a constant $c$ such that $\|f_i\|=c$ for all $i\in[n]$, then $\{f_i\}_{i\in[n]}$  is called a \emph{uniform frame}. If $\|f_i\|=1$ for all $i\in[n]$, then $\{f_i\}_{i\in[n]}$ is called a \emph{unit norm frame}.
	
	Let us see the definition of unitary equivalent frames.
	\begin{defn}
		Two frames $\{f_i\}_{i\in[n]}$ and $\{g_i\}_{i\in[n]}$ for $\mathbb{C}^k$ are said to be \emph{unitary equivalent} if there exists a unitary operator $U:\mathbb{C}^k\rightarrow\mathbb{C}^k$ such that $U(g_i)=f_i$ for all $i\in[n]$.
	\end{defn}
	Next, we see some operators associated with frames. If $\{f_i\}_{i\in[n]}$ is a Bessel sequence for the Hilbert space $\mathbb{C}^k$, then the \emph{analysis operator} $T:\mathbb{C}^k\rightarrow\mathbb{C}^n$ is defined by  
	\begin{align*}
		T(f)=\{\langle f,f_i\rangle\}_{i\in[n]}.
	\end{align*}
	The analysis operator is always linear and bounded. The Bessel sequence $\{f_i\}_{i\in[n]}$ is a frame if and only if the analysis operator $T$ is injective. The adjoint of $T$ is called the \emph{synthesis operator}, and it is defined by 
	\begin{align*}
		T^*(\{c_i\}_{i\in[n]})=\sum\limits_{i\in[n]}c_if_i.
	\end{align*}
	The canonical matrix representation $[T^*]$ of the synthesis operator $T^*$ of the frame $\{f_i\}_{i\in[n]}$ for $\mathbb{C}^k$ is a $k\times n$ matrix such that $[T^*]=[f_1\,\,f_2\,\,\cdots\,\,f_n]$.
	
	Associated with the frame $\{f_i\}_{i\in[n]}$, the \emph{frame operator} $S: \mathbb{C}^k\rightarrow \mathbb{C}^k$ is defined by 
	\begin{align*}
		S(f) = T^*T(f) =\sum_{i\in[n]}\langle f, f_i\rangle f_i.
	\end{align*}
	The frame operator $S$ is bounded,  linear, self-adjoint,  positive, and invertible. The frame operator provides the following reconstruction formula
	\begin{align*}
		f = SS^{-1}f \ =\sum_{i\in[n]}\langle f, S^{-1}f_i\rangle f_i=\sum_{i\in[n]}\langle f, f_i \rangle S^{-1}f_i,  \ f\in\mathbb{C}^k.
	\end{align*}
	Here, $\{S^{-1}f_i\}_{i\in[n]}$ is also a frame for $\mathbb{C}^k$, and it is called the \emph{canonical dual frame} of $\{f_i\}_{i\in[n]}$. If $n>k$, there exist infinitely many frames $\{g_i\}_{i\in[n]}$ that provide the reconstruction formula 
	\begin{align*}
		f =\sum_{i\in[n]}\langle f, f_i \rangle g_i \ =\sum_{i\in[n]}\langle f, g_i\rangle f_i,\ f\in\mathbb{C}^k.
	\end{align*}
	The frame $\{g_i\}_{i\in[n]}$ is called a \emph{dual frame} or an \emph{alternate dual frame} of $\{f_i\}_{i\in[n]}$. It is well known that $\{g_i\}_{i\in[n]}$ is a dual frame of $\{f_i\}_{i\in[n]}$ if and only if there exists a sequence of vectors $\{h_i\}_{i\in[n]}$ in $\mathbb{C}^k$ such that $\sum\limits_{i\in[n]}\langle f,h_i\rangle f_i=\textbf{0}$ for all $f\in\mathbb{C}^k$ and $g_i=S^{-1}f_i+h_i$ for all $i\in[n]$, where $S$ is the frame operator of $\{f_i\}_{i\in[n]}$.
	
	Another important operator associated with a frame is the \emph{Gramian operator}. If $\{f_i\}_{i\in[n]}$ is a frame for $\mathbb{C}^k$, then the Gramian (operator) $\mathcal{G}:\mathbb{C}^n\rightarrow\mathbb{C}^n$ is defined by $\mathcal{G}=TT^*$, where $T$ is the analysis operator of $\{f_i\}_{i\in[n]}$. The canonical matrix representation of the Gramian of a frame $\{f_i\}_{i\in[n]}$ is called the \emph{Gramian matrix} and is defined by 
	\begin{align*}
		\mathcal{G}=\left[\begin{array}{cccc}
			\|f_1\| ^2 & \langle f_2,f_1\rangle & \cdots &\langle f_n,f_1\rangle\\
			\langle f_1,f_2\rangle & \|f_2\|^2 &\cdots & \langle f_n,f_2\rangle\\
			\vdots & \vdots & \ddots & \vdots\\
			\langle f_1,f_n\rangle &  \langle f_2,f_n\rangle & \cdots &  \|f_n\|^2
		\end{array}\right].
	\end{align*}
	
	Next, we see the definition of strictly scalable frames.
	\begin{defn}\cite{CK1}
		A frame $\{f_i\}_{i\in[n]}$ for $\mathbb{C}^k$ is called \emph{strictly scalable} if there exist positive constants $\alpha_1,\alpha_2,\ldots,\alpha_n$ such that $\{\alpha_if_i\}_{i\in[n]}$ is a Parseval frame for $\mathbb{C}^k$.
	\end{defn}
	The following theorem gives the characterization of strictly scalable frames.
	\begin{thm}\cite{KOPT}\label{thm2.2}
		Suppose $\{f_i\}_{i\in[n]}$ is a frame for $\mathbb{C}^k$ with Gramian matrix $\mathcal{G}$. Then $\{f_i\}_{i\in[n]}$ is strictly scalable if and only if there exists an $n\times (n-k)$ matrix $M$ such that the matrix $\mathcal{G}+MM ^*$ is a positive definite diagonal matrix.
	\end{thm}
	For the general theory of frames and their applications, see \cite{BL, CK1, C2, DV, DV1}.

	\subsection{Matrices associated with graphs}
	Let $ G $ be a simple graph with vertex set $ V(G)=\{ v_1, v_2, \dots, v_n\}$. Two vertices $ v_i $ and $ v_j $ are adjacent if there is an edge between them, and we will write $ v_i \sim v_j $. Graph $G$ is a \emph{complete graph} if $v_i\sim v_j$ for all $i\neq j$. The degree of a vertex $ v_i $ is the number of edges adjacent to the vertex $ v_i $ and is denoted by $ d(v_i) $ or simply by $ d_i $. The smallest and the largest vertex degrees of $ G $ are denoted by $ \delta(G) $ (or simply by $\delta$) and $ \Delta(G) $ (or simply by $\Delta$), respectively. A vertex $v_i$ is called a \emph{null vertex} if $d(v_i)=0$. Graph $G$ is said to be a \emph{regular graph} if $d_i=d_j$ for all $i,j\in[n]$. If $d_i=d_j=r$ for all $i,j\in[n]$, then $G$ is called an \emph{$r$-regular graph}.  A $ p\times q $ matrix with all entries are ones is denoted by $ J_{p,q} $ (or simply by $ J $ when the order is clearly understood). The transpose and rank of a matrix $ M $ are denoted by $ M^t $ and $\text{rank}(M)$, respectively. If $M$ is a matrix of order $p\times q$, then $m_{ij}$ denotes the entry of the matrix $M$ in the $i^{\text{th}}$ row and $j^{\text{th}}$ column, and we write $M=[m_{ij}]_{1\leq i\leq p,1\leq j \leq q}$ (or simply $M=[m_{ij}]_{p\times q}$). Also, $diag(a_1,a_2,\ldots,a_n)$ is used to denote the $n\times n$ diagonal matrix with respective diagonal entries $a_1,a_2,\ldots,a_n$.
	
	Next, we see some matrices associated with a simple graph $G$. The \emph{degree matrix} of $G$, denoted by $D(G)$ (or simply by $D$), is an $n\times n$ diagonal matrix such that $D=diag(d_1,d_2,\ldots,d_n)$. The \emph{adjacency matrix} of the graph $G$, denoted by $A(G)=[a_{ij}]_{n\times n}$ (or simply by $A$), is an $n\times n$ matrix defined as
	\begin{align*}
		a_{ij}=\begin{cases}
			1, \ \text{if } i\neq j \text{ and } v_i\sim v_j\\
			0, \text{ otherwise.}
		\end{cases}
	\end{align*}
	The degree matrix and the adjacency matrix are symmetric matrices. By the eigenvalues of a graph $G$, we refer to the eigenvalues of the adjacency matrix $A(G)$. The \emph{Laplacian matrix} of a graph $G$, denoted by $L(G)$ (or simply by $L$), is an $n\times n$ symmetric matrix such that $L(G)=D(G)-A(G)$.  The Laplacian matrix is a positive semi-definite matrix. If $G$ is a simple graph of $n$ vertices with $p$ connected components, then the rank of $L(G)$ is $n-p$. 
	
	Next, we see a lower bound of the largest eigenvalues of $L(G)$.
	\begin{thm}\cite{Bapat}\label{thm2.3}
		Let $G$ be a graph with at least one edge and the largest vertex degree $\Delta$. If $\mu_1$ is the largest eigenvalues of the Laplacian matrix $L(G)$, then $\mu_1\geq \Delta+1$.
	\end{thm}
	
	It is well known that $0$ is always an eigenvalue of the Laplacian matrix, and it is the smallest eigenvalue. If $G$ is a graph with Laplacian matrix $L$ and $\mu_1\geq\mu_2\geq\cdots\geq\mu_n$ are the eigenvalues of $L$, then the second smallest eigenvalue, $\mu_{n-1}$, is called the \emph{algebraic connectivity} of $G$. Let us recall a lower bound of the algebraic connectivity.
	\begin{thm}\cite{Bapat}\label{thm2.4}
		Suppose $G$ is a connected graph of $n$ vertices with the smallest vertex degree $\delta$. If the eigenvalues of the Laplacian matrix $L(G)$ are $\mu_1\geq\mu_2\geq\ldots\geq\mu_n=0$, then $\mu_{n-1}\leq \frac{n}{n-1}\delta$.
	\end{thm}

	\section{Elementary properties of frames generated by graphs}
	
	We start this section with the definition of frames generated by a graph.
	
	\begin{defn}
		Suppose $\{f_i\}_{i\in[n]}$ is a frame for $\mathbb{C}^k$ with Gramian matrix $\mathcal{G}$. If $G$ is a simple graph with Laplacian matrix $L$ such that $\mathcal{G}=L$, then $\{f_i\}_{i\in[n]}$ is called a frame generated by graph $G$. In short, we call $\{f_i\}_{i\in[n]}$ as a \emph{$G(n,k)$-frame} for $\mathbb{C}^k$.
	\end{defn}
	
	Suppose $G$ is a simple graph with $n$ vertices and $n-k$ connected components. Then the rank of the Laplacian matrix $L$ is $k$. Now, we use the Laplacian matrix $L$ to generate a frame for the Hilbert space $\mathbb{C}^k$ as follows.
	
	\vspace{6pt}
	
	Since $L$ is a positive semi-definite matrix of rank $k$, thus $L$ has $k$ positive eigenvalues, and the remaining $n-k$ eigenvalues are $0$. Suppose the non-zero eigenvalues of $L$ are $\mu_1\geq\mu_2\geq\cdots\geq\mu_k$. Hence, there exist a real orthogonal matrix $M=[m_{ij}]_{1\leq i, j\leq n}$ and the diagonal matrix $D=\text{diag}(\mu_1,\mu_2,\ldots,\mu_k,0,\ldots,0)$ of order $n\times n$ such that 
	\begin{align*}
		L=M\, D\, M^{*}.
	\end{align*}
 
	Let us define the full column rank matrices $M_1$, $M_2$ and $D_1$ by  
	\begin{align*}
		M_1&=[m_{ij}]_{1\leq i \leq n, 1\leq j\leq k}\\ M_2&=[m_{ij}]_{1\leq i \leq n, k+1\leq j\leq n}\\ D_1&=\text{diag}\left(\sqrt{\mu_1},\ldots,\sqrt{\mu_k}\right).
	\end{align*}
	If $B=D_1M_1^*$, we have
	\begin{align*}
		L&=M\, D\, M^{*}\\
		&=\left[\begin{array}{cc}
			M_1   & M_2 
		\end{array}\right]\, \left[\begin{array}{cc}
			D_1^2   & \mathbf{0}\\
		\mathbf{0} & \mathbf{0}
		\end{array}\right]\, \left[\begin{array}{c}
			M_1^*\\
			M_2^*
		\end{array}\right]\\
		&=M_1\,D_1^2\,M_1^*\\
		&=B^*\,B.
	\end{align*}
	Since $M_1^*$ is a full row rank matrix and $D_1$ is invertible, $B$ is a full row rank matrix that is $\text{rank}(B)=k$.
	
	Let us generate a frame for $\mathbb{C}^k$ using the matrix $B$ as follows.
	
	\vspace{6pt}
	
	It is easy to observe that $\{B(e_i)\}_{i\in[n]}$ is a Bessel sequence for $\mathbb{C}^k$. Note that the matrix representation of the analysis operator $T:\mathbb{C}^k\rightarrow\mathbb{C}^n$ of the Bessel sequence $\{B(e_i)\}_{i\in[n]}$ is $B^*$. Since $B^*$ is a full column rank matrix so the analysis operator $T$ is injective, and hence the Bessel sequence $\{B(e_i)\}_{i\in[n]}$ is a frame for $\mathbb{C}^k$. Also, the Gramian matrix of the frame $\{B(e_i)\}_{i\in[n]}$ is $\mathcal{G}=[T]\,[T^*]=B^*B=L$.
	
	\begin{rem}
		If $G$ is a graph of $n$ vertices with Laplacian matrix $L$ of rank $k$, then the frame $\{B(e_i)\}_{i\in[n]}$ constructed above is a special kind of $G(n,k)$-frame and we call this frame as an \emph{$L_G(n,k)$-frame} for $\mathbb{C}^k$.
	\end{rem}

	Note that if $\{f_i\}_{i\in[n]}$ is a $G(n,k)$-frame for $\mathbb{C}^k$ with analysis operator $T$ such that $B$ is the canonical matrix representation of $T^*$, then the Laplacian matrix of the graph $G$ is $L=B^*B$ and $\{f_i\}_{i\in[n]}=\{B(e_i)\}_{i\in[n]}$.\\
	
	Next we see an example of an $L_G(4,3)$-frame for the Hilbert space $\mathbb{C}^3$.
	
	\begin{exa}
		Consider the graph $G$ given in FIGURE \ref{fig1}.
		
		\begin{figure}
			\begin{center}
				\includegraphics[scale= 0.65]{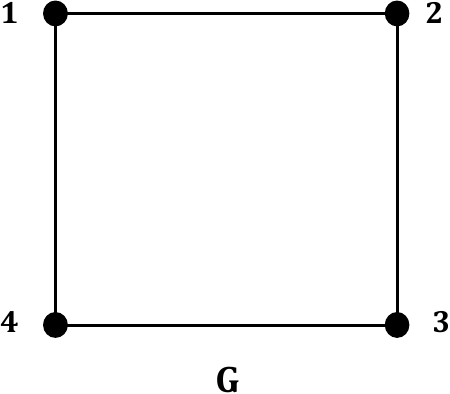}
				\caption{} \label{fig1}
			\end{center}
		\end{figure}
	
		The Laplacian matrix $L$ of $G$ is
		\begin{align*}
			L=\left[\begin{array}{cccc}
				2 &  -1 &0 &-1\\
				-1 & 2 & -1 & 0\\
				0 & -1 & 2 & -1\\
				-1 & 0 &-1 & 2
			\end{array}\right].
		\end{align*}
		The eigenvalues of $L$ are $4,\,2,\,2,$ and $0$. Next, we construct a real orthogonal matrix $M$ such that $L=M\,\text{diag}(4,2,2,0)\,M^*$. One may note that $[1\,\, -1\,\, 1\,\, -1]^t$, $[1\,\, 0\,\, -1\,\, 0]^t$, $[0\,\, 1\,\, 0\,\,-1]^t$ and $[1\,\,1\,\, 1\,\, 1]^t$ are the eigenvectors of $L$ corresponding to the eigenvalues $4,\,2,\,2$ and $0$, respectively. We can construct the matrix $M$ using the eigenvectors as follows
		\begin{align*}
			M=\left[\begin{array}{cccc}
				1/2 &  1/\sqrt{2} &0 & 1/2\\
				-1/2 & 0 & 1/\sqrt{2} & 1/2\\
				1/2 & -1/\sqrt{2} & 0 & 1/2\\
				-1/2 & 0 &-1/\sqrt{2} & 1/2
			\end{array}\right].
		\end{align*}
		It is easy to see that $M$ is a real orthogonal matrix. Also, $L=M\,\text{diag}(4,2,2,0)\,M^*$. Now choose the matrix $M_1$ as 
		\begin{align*}
			M_1=\left[\begin{array}{cccc}
				1/2 &  1/\sqrt{2} &0 \\
				-1/2 & 0 & 1/\sqrt{2} \\
				1/2 & -1/\sqrt{2} & 0\\
				-1/2 & 0 &-1/\sqrt{2} 
			\end{array}\right].
		\end{align*}
		Note that we form the matrix $M_1$ by taking the first three columns of the matrix $M$ that is $M_1$ contains the normalized eigenvectors of $L$ corresponding to the non-zero eigenvalues. Then
		\begin{align*}
			B&=\text{diag}(2,\sqrt{2},\sqrt{2}) M_1^*\\
			&=\left[\begin{array}{ccc}
				2 &  0 &0 \\
				0& \sqrt{2} & 0 \\
				0 & 0 & \sqrt{2}
			\end{array}\right]\left[\begin{array}{cccc}
				1/2 &  -1/2 & 1/2 & -1/2 \\
				1/\sqrt{2} & 0 & -1/\sqrt{2} & 0 \\
				0 & 1/\sqrt{2} & 0 & -1/\sqrt{2}
			\end{array}\right]\\
			&=\left[\begin{array}{cccc}
				1 &  -1 & 1 & -1 \\
				1 & 0 & -1 & 0 \\
				0 & 1 & 0 & -1
			\end{array}\right].
		\end{align*}
		Therefore, $\{B\,e_1,\,B\,e_2,\,B\,e_3,\,B\,e_4\}=\left\{\left[\begin{array}{c}
			1\\
			1\\
			0 
		\end{array}\right], \left[\begin{array}{c}
			-1\\
			0\\
			1 
		\end{array}\right], \left[\begin{array}{c}
			1\\
			-1\\
			0 
		\end{array}\right], \left[\begin{array}{c}
			-1\\
			0\\
			-1 
		\end{array}\right]\right\}$ is an $L_G(4,3)$-frame for $\mathbb{C}^3$.
	\end{exa}
	
	Every $L_G(n,k)$-frame is a $G(n,k)$-frame but the converse is not true. In the following example, we show a $G(n,k)$-frame which is not an $L_G(n,k)$-frame.
	
	\begin{exa}
		Consider the frame $\{f_i\}_{i\in[4]}=\left\{\left[\begin{array}{c}
			1\\
			1\\
			1 
		\end{array}\right],\, \left[\begin{array}{c}
			-1\\
			0\\
			0 
		\end{array}\right],\,\left[\begin{array}{c}
			0\\
			-1\\
			0
		\end{array}\right],\,\left[\begin{array}{c}
			0\\
			0\\
			-1 
		\end{array}\right]\right\}$ for the Hilbert space $\mathbb{C}^3$. Then the Gramian matrix of $\{f_i\}_{i\in[4]}$ is 
		\begin{align*}
			\mathcal{G}&=\left[\begin{array}{cccc}
				\|f_1\|^2 & \langle f_2,f_1\rangle & \langle f_3,f_1\rangle & \langle f_4,f_1\rangle\\
				\langle f_1,f_2\rangle & \|f_2\|^2  & \langle f_3,f_2\rangle & \langle f_4,f_2\rangle\\
				\langle f_1,f_3\rangle & \langle f_2,f_3\rangle & \|f_3\|^2  & \langle f_4,f_3\rangle\\
				\langle f_1,f_4\rangle & \langle f_2,f_4\rangle & \langle f_3,f_4\rangle & \|f_4\|^2 
			\end{array}\right]\\
			&=\left[\begin{array}{cccc}
				3 & -1 & -1 & -1\\
				-1 & 1  & 0 & 0\\
				-1 & 0 & 1 & 0\\
				-1 & 0 & 0 & 1 
			\end{array}\right].
		\end{align*}
		Consider the graph $G$ given in FIGURE \ref{fig2}.
		
		\begin{figure}
			\begin{center}
				\includegraphics[scale= 0.55]{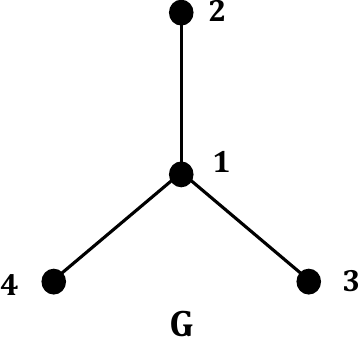}
				\caption{} \label{fig2}
			\end{center}
		\end{figure}
		
		The Laplacian matrix $L$ of $G$ is
		\begin{align*}
			L=\left[\begin{array}{cccc}
				3 & -1 & -1 & -1\\
				-1 & 1  & 0 & 0\\
				-1 & 0 & 1 & 0\\
				-1 & 0 & 0 & 1 
			\end{array}\right]=\mathcal{G}.
		\end{align*}
		Therefore, $\{f_i\}_{i\in[4]}$ is a $G(4,3)$-frame for $\mathbb{C}^3$. It is easy to observe that the columns of the canonical matrix representation of the analysis operator of an $L_G(n,k)$-frame are orthogonal to each other. The canonical matrix representation of the analysis operator of $G(4,3)$-frame $\{f_i\}_{i\in[4]}$ is 
		\begin{align*}
			\left[\begin{array}{ccc}
				1 & 1 & 1 \\
				-1 & 0 & 0 \\
				0 & -1 & 0 \\
				0 & 0 & -1 
			\end{array}\right].
		\end{align*}
		Since the columns of the above matrix are not orthogonal to each other, therefore $\left\{\left[\begin{array}{c}
			1\\
			1\\
			1 
		\end{array}\right],\, \left[\begin{array}{c}
			-1\\
			0\\
			0 
		\end{array}\right],\,\left[\begin{array}{c}
			0\\
			-1\\
			0
		\end{array}\right],\,\left[\begin{array}{c}
			0\\
			0\\
			-1 
		\end{array}\right]\right\}$ is a $G(4,3)$-frame for $\mathbb{C}^3$ but not an $L_G(4,3)$-frame.
	\end{exa}

	\begin{rem}
		For any positive integers $n>k$, there exists a $G(n,k)$-frame for $\mathbb{C}^k$. 
		
	For any positive integers $n>k$, we can always construct a simple graph $G$ with $n-k$ connected components. Then the rank of the Laplacian matrix $L$ of $G$ is $n-(n-k)=k$. Hence, the Laplacian matrix $L$ will generate an $L_G(n,k)$-frame for $\mathbb{C}^k$. 
	\end{rem}

	A graph $G$ can generate many frames for $\mathbb{C}^k$. In the following theorem, we show that $G(n,k)$-frames generated by the same graph $G$ are
	unitary equivalent.
	
	\begin{thm}
		If $G$ is a simple graph, then frames generated by $G$ are unitary equivalent.
	\end{thm}
	\proof
	Suppose $G$ is a simple graph with Laplacian matrix $L$. Let $F_1=\{f_i\}_{i\in[n]}$ and $F_2=\{g_i\}_{i\in[n]}$ be two frames generated by the graph $G$ with Gramian matrices $\mathcal{G}_1$ and $\mathcal{G}_2$, respectively. If the matrices $B=[f_1\,\,f_2\,\,\cdots\,\,f_n]$ and $C=[g_1\,\,g_2\,\,\cdots\,\,g_n]$, then $F_1=\{B(e_i)\}_{i\in[n]}$ and $F_2=\{C(e_i)\}_{i\in[n]}$. Thus, $\mathcal{G}_1=B^*B=L=\mathcal{G}_2=C^*C$ that is $B^*B=C^*C$. Let $S_1$ and $S_2$ be the frame operators of $F_1$ and $F_2$, respectively. Then $(BC^*)(CB^*)=B(C^*C)B^*=(BB^*)(BB^*)=S_1^2$ which gives $S_1^{-2}=(CB^*)^{-1}(BC^*)^{-1}$. Take $U=S_1^{-2}BC^*S_2$. Then we have
	\begin{align*}
		U^*U
		&=S_2CB^*S_1^{-2}S_1^{-2}BC^*S_2\\
		&=S_2CB^*(CB^*)^{-1}(BC^*)^{-1}(CB^*)^{-1}(BC^*)^{-1}BC^*S_2\\
		&=S_2(CB^*BC^*)^{-1}S_2\\
		&=S_2(CC^*CC^*)^{-1}S_2\\
		&=S_2S_2^{-2}S_2=I.
	\end{align*}
	Similarly, we have
	\begin{align*}
		UU^*&=S_1^{-2}BC^*S_2S_2CB^*S_1^{-2}\\
		&=S_1^{-2}BC^*CC^*CC^*CB^*S_1^{-2}\\
		&=S_1^{-2}BB^*BB^*BB^*BB^*S_1^{-2}\\
		&=S_1^{-2}S_1^4S_1^{-2}=I.
	\end{align*}
	Therefore, $U$ is a unitary operator. Further, we have
	\begin{align*}
		UC&=S_1^{-2}BC^*S_2C\\
		&=S_1^{-2}BC^*CC^*C\\
		&=S_1^{-2}BB^*BB^*B\\
		&=S_1^{-2}S_1^2B=B.\\
	\end{align*}
	Thus, $Ug_i=UC(e_i)=B(e_i)=f_i$ for all $i\in[n]$. Hence, the frames $F_1$ and $F_2$ are unitary equivalent.
	\endproof
	
	In the next proposition, we show that unitary equivalent $G(n,k)$-frames must be generated by the same graph.
	
	\begin{prop}
		Suppose $F_1$ and $F_2$ are $G_1(n,k)$-frame and $G_2(n,k)$-frame, respectively. If the frames $F_1$ and $F_2$ are unitary equivalent, then $G_1$ and $G_2$ are identical graphs.
	\end{prop}
	\proof
	Suppose $T_1$ and $T_2$ are the analysis operators of the frames $F_1$ and $F_2$, respectively. Since $F_1$ and $F_2$ are unitary equivalent, therefore there exists a unitary operator $U:\mathbb{C}^k\rightarrow\mathbb{C}^k$ such that $T_1=T_2U$. Let $L(G_1)$ and $L(G_2)$ be the Laplacian matrices of $G_1$ and $G_2$, respectively. If $\mathcal{G}_1$ and $\mathcal{G}_2$ are the Gramian matrices of the frames $F_1$ and $F_2$, respectively, then
	\begin{align*}
		L(G_1)=\mathcal{G}_1=T_1T_1^*=(T_2U)(T_2U)^*=T_2UU^*T_2^*=T_2T_2^*=\mathcal{G}_2=L(G_2).
	\end{align*}
	Thus, $L(G_1)=L(G_2)$. Therefore, $G_1$ and $G_2$ are identical graphs.
	\endproof 
	It is known that the Laplacian matrices of isomorphic graphs need not be the same. Thus, frames generated by isomorphic graphs may not be unitary equivalent.
	\begin{rem}
		Distinct (non-identical) graphs generate distinct unitary equivalence classes of frames.
	\end{rem}
	
	Canonical dual and alternate dual frames are used in the reconstruction of vectors or signals from frame coefficients. Next, we study about the canonical dual and alternate dual frames of an $L_G(n,k)$-frame and $G(n,k)$-frame.  First, we prove the following lemma which provides the frame operator of an $L_G(n,k)$-frame.
	
	\begin{lem}\label{lem1.9}
		If $F$ is an $L_G(n,k)$-frame, then the frame operator $S$ of the frame $F$ is a diagonal matrix. Further, if $\mu_1\geq\mu_2\geq\cdots\geq\mu_k$ are the non-zero eigenvalues of the Laplacian matrix of $G$, then $S=diag(\mu_1,\mu_2,\ldots,\mu_k)$.
	\end{lem} 
	\proof
	It is given that $F$ is an $L_G(n,k)$-frame, thus the rank of the Laplacian matrix $L$ of $G$ is $k$. Suppose the non-zero eigenvalues of $L$ are $\mu_1\geq\mu_2\geq\cdots\geq\mu_k$. Since $F$ is an $L_G(n,k)$-frame, thus $F=\{B(e_i)\}_{i\in[n]}$ where $B=diag(\sqrt{\mu_1},\sqrt{\mu_2},\ldots,\sqrt{\mu_k})M_1^*$ and $L_G=M\, diag(\mu_1,\mu_2,\ldots,\mu_k,0,\ldots,0)\,M^*$ such that $M$ is a real orthogonal matrix and $M_1$ is obtained by taking the first $k$ columns of $M$. Since $M$ is orthogonal so $M_1^*M_1=I$. 
	
	The frame operator $S$ of the frame $F$ is given by $S(f)=\sum\limits_{i\in[n]}\langle f,B(e_i)\rangle B(e_i)=BB^*(f)$. Thus, $S=BB^*$. Then we have
	\begin{align*}
		S&=diag(\sqrt{\mu_1},\sqrt{\mu_2},\ldots,\sqrt{\mu_k})M_1^*M_1 diag(\sqrt{\mu_1},\sqrt{\mu_2},\ldots,\sqrt{\mu_k})\\
		&=diag(\sqrt{\mu_1},\sqrt{\mu_2},\ldots,\sqrt{\mu_k})diag(\sqrt{\mu_1},\sqrt{\mu_2},\ldots,\sqrt{\mu_k})\\
		&=diag(\mu_1,\mu_2,\ldots,\mu_k).
	\end{align*}
	Therefore, $S$ is a diagonal matrix.
	\endproof
	
	In the next theorem, we obtain the canonical dual frame of an $L_G(n,k)$-frame. 
	\begin{thm}\label{thm1.10}
		Suppose $F=\{f_i\}_{i\in[n]}$ is an $L_G(n,k)$-frame generated by the graph $G$ and $\mu_1\geq\mu_2\geq\cdots\geq\mu_k$ are the non-zero eigenvalues of the Laplacian matrix of $G$. Then canonical dual frame of $F$ is $\{Df_1,Df_2,\ldots,Df_n\}$ where $D=diag(\mu_1^{-1},\mu_2^{-1},\ldots,\mu_k^{-1})$.
	\end{thm}
	\proof
	Let $S$ be the frame operator of $F$. Then by Lemma \ref{lem1.9}, $S=diag(\mu_1,\mu_2,\ldots,\mu_k)$. Thus, $S^{-1}=diag(\mu_1^{-1},\mu_2^{-1},\ldots,\mu_k^{-1})=D$. Hence, the canonical dual frame of $F$ is $$S^{-1}F=\{S^{-1}f_1,S^{-1}f_2,\ldots,S^{-1}f_n\}=\{Df_1,Df_2,\ldots,Df_n\}.$$
	\endproof
	
	For $n>k$, there exist infinitely many dual frames of a $G(n,k)$-frame. In the following theorem, we obtain the family of alternate dual frames of the frames generated by connected graphs. Note that if $G$ is a connected graph with $n$ vertices, then the rank of its Laplacian matrix is $n-1$, and hence we get a $G(n,n-1)$-frame for $\mathbb{C}^{n-1}$.

	\begin{thm}\label{thm1.11}
		Suppose $G$ is a connected graph with $n$ vertices and $F=\{f_i\}_{i\in[n]}$ is a $G(n,n-1)$-frame for $\mathbb{C}^{n-1}$ with frame operator $S$. Then the family of dual frames of $F$ is $\left\{\{S^{-1}f_i+h\}_{i\in[n]}:h\in\mathbb{C}^{n-1}\right\}$.
	\end{thm}
	\proof
	It is known that the dual frames of $F$ are of the form $\{S^{-1}f_i+h_i\}_{i\in[n]}$ where $\{h_i\}_{i\in[n]}\subset\mathbb{C}^{n-1}$ satisfy $\sum\limits_{i\in[n]}\langle f,h_i\rangle f_i=\textbf{0}$ for all $f\in\mathbb{C}^{n-1}$. Thus, we have
	\begin{align}\label{eq1.2}
		[f_1\ f_2\ \cdots\ f_n]\left[\begin{array}{c}
			\langle f,h_1\rangle \\
			\langle f,h_2\rangle  \\
			\vdots \\
			\langle f,h_n\rangle  
		\end{array}\right]=\textbf{0},\text{ for all } f\in\mathbb{C}^{n-1}.
	\end{align}
	Suppose $T$ is the analysis operator of the frame $F$. Then $f_i=T^*(e_i)$ for all $i\in[n]$ and the Gramian matrix $\mathcal{G}=L=TT^*$, where $L$ is the Laplacian matrix of the graph $G$. Thus, by Equation \eqref{eq1.2}, $T^*h_f=\textbf{0}$ where $h_f=\left[\begin{array}{c}
		\langle f,h_1\rangle \\
		\langle f,h_2\rangle  \\
		\vdots \\
		\langle f,h_n\rangle  
	\end{array}\right]$. This gives that $Lh_f=\textbf{0}$. Since $G$ is a connected graph, thus the rank of $L$ is $n-1$. Hence, the nullity of $L$ is $n-(n-1)=1$. Since $L\left[\begin{array}{cccc}
		1 & 1 & \cdots & 1
	\end{array}\right]^t=\textbf{0}$,  we have $$
	\text{null space of } L=span\left\{\left[\begin{array}{c}
		1 \\
		1 \\
		\vdots \\
		1 
	\end{array}\right]\right\}.$$ Hence, $h_f=\alpha \left[\begin{array}{cccc}
		1 & 1 & \cdots & 1
	\end{array}\right]^t$ for some scalar $\alpha$. This gives $\langle f,h_i\rangle=\langle f,h_j\rangle$ for all $i, j\in[n]$ and for all $f\in\mathbb{C}^{n-1}$. Hence, $h_i=h_j$ for all $i,j\in[n]$. Therefore, the family of dual frames of $F$ is $\left\{\{S^{-1}f_i+h\}_{i\in[n]}:h\in\mathbb{C}^{n-1}\right\}$.
	\endproof
	
	In the following corollary, we obtain the dual frame of an $L_G(n,k)$-frame generated by the connected graph $G$.
	\begin{cor}
		Suppose $G$ is a connected graph with $n$ vertices and $F$ is an $L_G(n,n-1)$-frame. If $\mu_1\geq\mu_2\geq\cdots\geq\mu_{n-1}>\mu_n=0$ are the eigenvalues of the Laplacian matrix of $G$ and $D=diag(\mu_1^{-1},\mu_2^{-1},\ldots,\mu_{n-1}^{-1})$, then the family of dual frames of $F$ is $\{\{Df_i+h\}_{i\in[n]}:h\in\mathbb{C}^{n-1}\}$.
	\end{cor}
	\proof
	By using Theorem \ref{thm1.10} and Theorem \ref{thm1.11}, we have the family of dual frames of $F$ is $\{\{Df_i+h\}_{i\in[n]}:h\in\mathbb{C}^{n-1}\}$.
	\endproof
	
	In the next theorem, we obtain the family of dual frames of a $G(n,k)$-frame where $G$ is any graph. Here, the graph $G$ need not be a connected graph. 
	
	\begin{thm}\label{thm1.13}
		Suppose $G$ is a simple graph with connected components $G_1,G_2,\ldots,\break G_{m}$ such that $G_1,G_2,\ldots,G_{m}$ has vertices $\{v_1,v_2,\ldots,v_{n_1}\}, \{v_{n_1+1}, v_{n_1+2}, \ldots, v_{n_2}\}, \break\ldots,\{v_{n_{m-1}+1}, v_{n_{m-1}+2}, \ldots, v_{n_m}=v_n\}$, respectively, where $m=n-k$. If $F=\{f_i\}_{i\in[n]}$ is a $G(n,k)$-frame with the frame operator $S$, then any dual frame of $F$ is of the form $\{S^{-1}f_i+\nu_1\}_{i=1}^{n_1}\bigcup\{S^{-1}f_i+\nu_2\}_{i=n_1+1}^{n_2}\bigcup\cdots\bigcup\{S^{-1}f_i+\nu_m\}_{i=n_{m-1}+1}^{n}$ where $\nu_1,\nu_2,\ldots,\nu_m$ are arbitrary vectors in $\mathbb{C}^k$.
	\end{thm}
	\proof
	It is known that any dual frame of $F$ is of the form $\{S^{-1}f_i+h_i\}_{i\in[n]}$ such that $\{h_i\}_{i\in[n]}\subset\mathbb{C}^k$ satisfy $\sum\limits_{i\in[n]}\langle f,h_i\rangle f_i=\textbf{0}$. Since the connected components of the graph $G$ are $G_1,G_2,\ldots,G_m$, thus the Laplacian matrix $L(G)$ of the graph $G$ is 
	\begin{align*}
		L(G)=\left[\begin{array}{cccc}
			L(G_1) & \mathbf{0} & \cdots & \mathbf{0}\\
			\mathbf{0} & L(G_2) & \cdots & \mathbf{0}\\
			\vdots & \vdots & \ddots & \vdots\\
			\mathbf{0} & \mathbf{0} & \cdots & L(G_m)
		\end{array}\right].
	\end{align*}
 Since $\sum\limits_{i\in[n]}\langle f,h_i\rangle f_i=\textbf{0}$, we have  $$L(G)\left[\begin{array}{c}
		\langle f,h_1\rangle\\
		\langle f,h_2\rangle\\
		\vdots\\
		\langle f,h_n\rangle
	\end{array}\right]=\textbf{0}.$$ 
	Thus, for $i\in[m]$, $L(G_i)\left[\begin{array}{c}
		\langle f,h_{n_{i-1}+1}\rangle\\
		\langle f,h_{n_{i-1}+2}\rangle\\
		\vdots\\
		\langle f,h_{n_{i}}\rangle
	\end{array}\right]=\textbf{0}$, where $n_0=0$. Since, $G_i$ is connected, nullity of  $L(G_i)$ is $1$ and the null space of $L(G_i)$ is $span\{[1\, 1\, \cdots\, 1]^t\}$. Hence, $\langle f, h_{n_{i-1}+1}\rangle=\langle f, h_{n_{i-1}+2}\rangle=\cdots=\langle f, h_{n_{i}}\rangle$ for all $f\in\mathbb{C}^k$. Then, for any $i\in[m]$, we have
	\begin{align*}
		h_{n_{i-1}+1}=h_{n_{i-1}+2}=\cdots=h_{n_{i}}.
	\end{align*}
	Therefore, any dual frame of $F$ is of the form $\{S^{-1}f_i+\nu_1\}_{i=1}^{n_1}\cup\{S^{-1}f_i+\nu_2\}_{i=n_1+1}^{n_2}\cup\cdots\cup\{S^{-1}f_i+\nu_m\}_{i=n_{m-1}+1}^{n}$, where $\nu_1,\nu_2,\ldots,\nu_m$ are arbitrary vectors in $\mathbb{C}^k$.
	\endproof

	\section{Tight frames and regular graphs}
	
	In signal processing, the reconstruction of a signal is numerically optimal stable from tight frame coefficients \cite{CK1}. Tight frames are also known to be optimal for handling erasures. In this section, we characterize tight frames generated by graphs in terms of regular graphs. Regular graphs are one of the important classes of graphs. They are widely studied in graph theory due to many aspects, including studying Ramanujan graphs and their connection with the Riemann hypothesis, see \cite{MSS1}.  
	
	In the following theorem, we present a necessary condition of tight frames generated by graphs.
	
	\begin{thm}\label{thm4.1}
		Suppose $G$ is a graph of $n$ vertices and $F=\{f_i\}_{i\in[n]}$ is a $G(n,k)$-frame. If $F$ is a tight frame for $\mathbb{C}^k$, then each connected component of $G$ is regular.
	\end{thm}
	\proof
	Suppose $F$ is an $\alpha$-tight frame and the graph $G$ has connected components $G_1,G_2,\ldots,G_{n-k}$. Let the number of vertices in the component $G_i$ be $n_i$ for $i\in[n-k]$. The rank of $L(G)$ is $k\, (=n-(n-k))$. Thus, the frame operator of the frame $F$ is $S=\alpha I$, where $I$ is the identity map on $\mathbb{C}^k$.
	
	The Gramian matrix $\mathcal{G}$ of the frame $F$ is given by  
	\begin{align*}
		\mathcal{G}=L(G)=\left[\begin{array}{cccc}
			L(G_1) &  \mathbf{0} & \cdots & \mathbf{0} \\
			\mathbf{0} & L(G_2) & \cdots & \mathbf{0} \\
			\vdots & \vdots & \ddots & \vdots\\
			\mathbf{0} & \mathbf{0} & \cdots & L(G_{n-k})
		\end{array}\right].
	\end{align*}
	Since the non-zero eigenvalues of the Gramian matrix $\mathcal{G}$ are the same as the non-zero eigenvalues of the frame operator $S$, thus the eigenvalues of $L(G)$ are $\alpha$ with multiplicity $k$ and $0$ with multiplicity $n-k$.
	
	As $G_i$ is connected, $L(G_i)$ has $n_i-1$ non-zero eigenvalues and one $0$ eigenvalue for $i\in[n-k]$. If $\mu_1^{(i)}\geq\mu_2^{(i)}\geq\ldots\geq\mu_{n_i}^{(i)}$ are the eigenvalues of $L(G_i)$, then $\mu_1^{(i)}=\mu_2^{(i)}=\cdots=\mu_{n_i-1}^{(i)}=\alpha$ and $\mu_{n_i}^{(i)}=0$. If $\Delta_i$ and $\delta_i$ are the largest and smallest vertex degrees of $G_i$, respectively, then by the Theorem \ref{thm2.3} and Theorem \ref{thm2.4}, $\mu_1^{(i)}=\alpha\geq\Delta_i+1$ and $\mu_{n_i-1}^{(i)}=\alpha\leq \frac{n_i}{n_i-1}\delta_i$. That is 
	\begin{align*}
		\Delta_i+1\leq\alpha\leq\frac{n_i}{n_i-1}\delta_i.
	\end{align*}
	Thus, $n_i\Delta_i+(n_i-1-\Delta_i)\leq n_i\delta_i$. Since $G$ is a simple graph, $\Delta_i\leq n_i-1$. Hence, $n_i\Delta_i\leq n_i\delta_i$ that is $\Delta_i\leq\delta_i$. Then $\Delta_i=\delta_i$. Hence, $G_i$ is regular. Therefore, each connected component of the graph $G$ is regular.
	\endproof
	
	In the next corollary, we show that tight frames cannot be generated by non-regular connected graphs. 
	\begin{cor}\label{cor0.7}
		Suppose $G$ is a connected graph of $n$ vertices and $F=\{f_i\}_{i\in[n]}$ is a $G(n,n-1)$-frame for $\mathbb{C}^{n-1}$. If $F$ is a tight frame, then $G$ must be a regular graph.
	\end{cor} 
	\proof
	Since $F$ is a tight frame, by the Theorem \ref{thm4.1}, each connected component of $G$ is regular. It is given that $G$ is a connected graph and hence $G$ is a regular graph.
	\endproof
	
	In Theorem \ref{thm4.1}, we show that if a frame generated by a graph $G$ is a tight frame, then each connected component of $G$ is regular. In the next theorem, we prove a stronger result by showing that if a tight frame is generated by a graph with no null vertex, then the graph must be a regular graph. 
	
	\begin{thm}\label{thm4.3}
		Suppose $G$ is a graph of $n$ vertices with no null vertex and $F=\{f_i\}_{i\in[n]}$ is a $G(n,k)$-frame for $\mathbb{C}^k$. If $F$ is a tight frame, then $G$ is a regular graph.
	\end{thm}
	\proof
	Suppose $F$ is an $\alpha$-tight frame and the graph $G$ has connected components $G_1,G_2,\ldots,G_{n-k}$. Then by the Theorem \ref{thm4.1}, $G_i$ is regular for all $i\in[n-k]$. 
	
	Suppose $G_i$ is an $r_i$-regular graph for $i\in[n-k]$. Since $G$ has no null vertex, $r_i>0$ for all $i\in[n-k]$. The adjacency matrix $A(G)$ of $G$ is 
	\begin{align*}
		A(G)=\left[\begin{array}{cccc}
			A(G_1) &  \mathbf{0} & \cdots & \mathbf{0} \\
			\mathbf{0} & A(G_2) & \cdots & \mathbf{0} \\
			\vdots & \vdots & \ddots & \vdots\\
			\mathbf{0} & \mathbf{0} & \cdots & A(G_{n-k})
		\end{array}\right]
	\end{align*}
	where $A(G_i)$ is the adjacency matrix of $G_i$.  Then
	the Gramian matrix $\mathcal{G}$ of the frame $F$ is
	\begin{align*}
		\mathcal{G}=L(G)=D(G)-A(G)=\left[\begin{array}{cccc}
			r_1I-A(G_1) &  \mathbf{0} & \cdots & \mathbf{0} \\
			\mathbf{0} & r_2I-A(G_2) & \cdots & \mathbf{0} \\
			\vdots & \vdots & \ddots & \vdots\\
			\mathbf{0} & \mathbf{0} & \cdots & r_kI-A(G_{n-k})
		\end{array}\right]
	\end{align*}
	where $I$ is an identity matrix of appropriate order.
	
	If $S$ is the frame operator and $T$ is the analysis operator of the frame $F$, then $S=T^*T=\alpha I$. Thus, $\mathcal{G}^2=TT^*TT^*=TST^*=T\alpha IT^*=\alpha \mathcal{G}$ that is $\mathcal{G}^2=\alpha \mathcal{G}$. Then we have
	\begin{align*}
		&\left[\begin{array}{cccc}
			(r_1I-A(G_1))^2 &  \mathbf{0} & \cdots & \mathbf{0} \\
			\mathbf{0} & (r_2I-A(G_2))^2 & \cdots & \mathbf{0} \\
			\vdots & \vdots & \ddots & \vdots\\
			\mathbf{0} & \mathbf{0} & \cdots & (r_kI-A(G_{n-k}))^2
		\end{array}\right]\\
		&\quad \quad=\left[\begin{array}{cccc}
			\alpha(r_1I-A(G_1)) &  \mathbf{0} & \cdots & \mathbf{0} \\
			\mathbf{0} & \alpha(r_2I-A(G_2)) & \cdots & \mathbf{0} \\
			\vdots & \vdots & \ddots & \vdots\\
			\mathbf{0} & \mathbf{0} & \cdots & \alpha(r_kI-A(G_{n-k}))
		\end{array}\right].
	\end{align*}
	Thus, $(r_i I-A(G_i))^2=\alpha (r_i I-A(G_i))$ that is $r_i^2 I+A(G_i)^2-2 r_i A(G_i)=\alpha r_i I-\alpha A(G_i)$. Then
	\begin{align}\label{eq1}
		A(G_i)^2=(\alpha r_i-r_i^2)I+(2 r_i-\alpha)A(G_i).
	\end{align}
	Since the diagonal entries of the matrix $A(G_i)^2$ are the degrees of the corresponding vertices of $G_i$ and $G_i$ is an $r_i$-regular graph, thus $r_i=\alpha r_i-r_i^2$. This gives that $\alpha=r_i+1$. Thus, $G_i$ is an $(\alpha-1)$-regular graph for all $i\in[n-k]$. Therefore, $G$ is an $(\alpha-1)$-regular graph.
	\endproof
	
	Uniform tight frames are used to handle erasures, see \cite{BP, CK}. In general, it is not easy to find uniform tight frames. In the following corollary, we show that tight frames generated by graphs with no null vertex are always uniform frames.
	
	\begin{cor}
		Suppose $G$ is a graph of $n$ vertices with no null vertex and $F=\{f_i\}_{i\in[n]}$ is a $G(n,k)$-frame. If $F$ is a tight frame, then $F$ must be a uniform frame.
	\end{cor}
	\proof
	Since $G$ is a graph with no null vertex and $F$ is a tight frame, therefore by the Theorem \ref{thm4.3}, $G$ is a regular graph. Let $G$ be an $r$-regular graph. Then the Gramian matrix $\mathcal{G}=[g_{ij}]_{n\times n}=L=r I-A$, where $A$ and $L$ are the adjacency and Laplacian matrices of $G$, respectively. Thus, $\|f_i\|^2=g_{ii}=r$ for all $i\in[n]$. Hence, $F$ is a uniform frame.
	\endproof
	
	Next, we provide a characterization of tight $G(n,k)$-frames in terms of graphs with exactly two distinct eigenvalues.
	
	\begin{thm}\label{thm4.5}
		Suppose $G$ is a graph with no null vertex and $F=\{f_i\}_{i\in[n]}$ is a $G(n,k)$-frame for $\mathbb{C}^k$. Then $F$ is a tight frame if and only if $G$ has exactly two distinct eigenvalues.
	\end{thm}
	\proof
	Let $G$ be the graph with $n$ vertices and connected components $G_1,\ldots,G_{n-k}$. 
	
	First, assume that $F$ is an $\alpha$-tight frame. Then by the Theorem \ref{thm4.3}, $G$ is a regular graph, say $r$-regular graph. If $A$ is the adjacency matrix of $G$, then the Gramian matrix of the frame $F$ is $\mathcal{G}=rI-A$. Since $F$ is an $\alpha$-tight frame, the frame operator $S=\alpha I$ and the Gramian matrix satisfies $\mathcal{G}^2=\alpha \mathcal{G}$. Thus, we have $ (rI-A)^ 2=\alpha(rI-A) $. That is $ A^2-(2r-\alpha)A+r(r-\alpha)I=\textbf{0}$. Thus, $A^2-(2r-\alpha)A+r(r-\alpha)I=\textbf{0}$ is the minimal polynomial of $A$. Hence, the eigenvalues of $A$ are $r$ and $r-\alpha$. A simple calculation shows that the multiplicity of the eigenvalue $r$ is $\frac{n(\alpha-r)}{\alpha}$ and multiplicity of the eigenvalue $r-\alpha$ is $\frac{nr}{\alpha}$. Hence, $G$ is a graph with exactly two distinct eigenvalues. 
	
	Conversely, suppose that $G$ has exactly two distinct eigenvalues $\alpha$ and $\beta$ (say). Then $G$ must be a regular graph, say $r$-regular graph. Since $r$ must be an eigenvalue of $G$, let $\beta=r$. 
	
	If $A$ is the Adjacency matrix of $G$, then the minimal polynomial of $A$ is $(A-rI)(A-\alpha I)=\textbf{0}$ that is $A^2-(\alpha+r)A+r\alpha I=\textbf{0}$. Thus, $A^2=(\alpha+r)A-r\alpha I$. The Gramian matrix $\mathcal{G}=rI-A$. Then we have
	\begin{align*}
		\mathcal{G}^2&=(rI-A)^2\\
		&=A^2-2rA+r^2I\\
		&=(\alpha+r)A-r\alpha I-2rA+r^2I\\
		&=(\alpha-r)A-r(\alpha-r) I\\
		&=(r-\alpha)\mathcal{G}.
	\end{align*}
	Let $T$ and $S$ be the analysis operator and frame operator of the frame $F$, respectively. Since $F$ is a frame, therefore the  analysis operator $T$ is injective and hence it has a left inverse. Then we have
	\begin{align*}
		S&=T^*T\\
		&=T^{-1}TT^*TT^*(T^*)^{-1}\\
		&=T^{-1}\mathcal{G}^2(T^*)^{-1}\\
		&=T^{-1}(r-\alpha)\mathcal{G}(T^*)^{-1}\\
		&=(r-\alpha)I.
	\end{align*}
	Therefore, $F$ is an $(r-\alpha)$-tight frame for $\mathbb{C}^k$.
	\endproof

	Many authors studied frames using graph-theoretic properties, see \cite{BP, SH}. In the following theorem, we prove the well-known graph property (connected graphs with exactly two distinct eigenvalues are complete graphs) using the concept of scalable frames.
	\begin{thm}\label{thm4.6}
		Suppose $G$ is a connected graph. Then $G$ has only two distinct eigenvalues if and only if $G$ is a complete graph.
	\end{thm}
	\proof
	First, assume that $G$ is a connected graph of $n$ vertices with exactly two distinct eigenvalues. Suppose that $L,\, A$, and $D$ are the Laplacian matrix, Adjacency matrix, and degree matrix of the graph $G$, respectively. Let $F=\{f_i\}_{i\in[n]}$ be an $L_G(n,n-1)$-frame for $\mathbb{C}^{n-1}$. By the Theorem \ref{thm4.5}, $F$ is a tight frame for $\mathbb{C}^{n-1}$.
	
	Suppose that $F$ is an $\alpha$-tight frame. Then $\frac{1}{\sqrt{\alpha}}F=\left\{\frac{1}{\sqrt{\alpha}}f_i\right\}_{i\in[n]}$ is a Parseval frame for $\mathbb{C}^{n-1}$. Hence, $F$ is a strictly scalable frame. Thus, by the Theorem \ref{thm2.2}, there exists a matrix $R$ of order $n\times 1$ such that $\mathcal{G}+RR^*$ is a positive definite diagonal matrix, where $\mathcal{G}$ is the Gramian matrix of the frame $F$. Let $D_1=\mathcal{G}+RR^*=L+RR^*=D-A+RR^*$. Then $D_1$ is a positive definite diagonal matrix. Thus, $RR^*=(D_1-D)+A=D_2+A$ where $D_2=D_1-D$ is a diagonal matrix. We know $rank(RR^*)=rank(R)\leq 1$. Since $G$ is a connected graph, $D_2+A\neq \textbf{0}$. Hence, $rank(D_2+A)=rank(RR^*)\geq 1$. Thus, $rank(D_2+A)=rank(RR^*)=1$. Since $G$ is a connected graph, each row of the matrix $D_2+A$ must be non-zero. If $r_1$ is the first row of the matrix $D_2+A$, then there exist scalars $\alpha_1,\ldots,\alpha_{n-1}$ such that
	$$D_2+A=\left[\begin{array}{c}
		r_1\\
		\alpha_1 r_1\\
		\vdots\\
		\alpha_{n-1}r_1
	\end{array}\right].$$
	Since $A$ is a symmetric matrix, $D_2+A$ is also a symmetric matrix. Thus, $\alpha_1=\cdots=\alpha_{n-1}=1$. Hence, we have $$D_2+A=\left[\begin{array}{c}
		r_1\\
		r_1 \\
		\vdots\\
		r_1
	\end{array}\right].$$
	If there exist $i\neq j$ such that $a_{ij}=0$, then each entry of the $j^{th}$column of $D_2+A$ is $0$. Thus, $a_{ij}=0$ for all $i\in[n]$ and hence $G$ is a disconnected graph which is not possible. Thus, every non-diagonal entry of the adjacency matrix $A$ must be $1$. Therefore, $G$ is a complete graph.

	Conversely, assume that $G$ is a complete graph. Then $A=J-I$. Thus, $A^2=(J-I)^2=J^2+I-2J=nJ-2J+I$. Hence, the minimal polynomial of $A$ is $A^2-(n-2)A-(n-1)I=\textbf{0}$. Therefore, $G$ has only two distinct eigenvalues.
	
	\endproof
	
	We obtain the following characterization of tight frames generated by connected graphs using Theorem \ref{thm4.5} and Theorem \ref{thm4.6}. The characterization shows that tight frames can be generated using complete graphs.
	\begin{cor}\label{cor2.6}
		Suppose $G$ is a connected graph and $F$ is a frame generated by $G$. Then $F$ is a tight frame if and only if $G$ is a complete graph.
	\end{cor}

	\mbox{}
	
	\endproof

\begin{thebibliography}{99}\baselineskip12pt
		
		\bibitem{NAS}
		F. Arabyani-Neyshaburi, A. A. Arefijamaal and G. Sadeghi, Numerically and spectrally optimal dual frames in Hilbert spaces, \emph{Linear Algebra Appl.}, 604 (2020), 52--71.
		
		\bibitem{AO}
		T. Adali and A. Ortega, Applications of Graph Theory [Scanning the Issue], \emph{Proceedings of the IEEE}, 106 (5) (2018), 784--786.
		
		\bibitem{Bapat}
		R. B. Bapat, Graphs and matrices, Springer, London, Hindustan Book Agency, New Delhi, 2010.
		
		\bibitem{BL}
		J. Benedetto and S. Li, The theory of multiresolution analysis frames and applications to filter banks, \emph{Appl. Comput. Harmon. Anal.}, 5 (1998), 389--427.
		
		\bibitem{B}
		N. Biggs, Algebraic graph theory, Second edition, Cambridge University Press, Cambridge, 1993.
		
		\bibitem{BP}
		B. G. Bodmann and V. I. Paulsen, Frames, graphs and erasures, \emph{Linear Algebra Appl.}, 404 (2005), 118--146.
		
		\bibitem{BP1}
		B. G. Bodmann and V. I. Paulsen, Frame paths and error bounds for sigma–delta quantization, \emph{Appl. Comput. Harmon. Anal.}, 22 (2) (2007), 176--197.
		
		
		\bibitem{BH}
		A. E. Brouwer and W. H.  Haemers, Spectra of graphs, Universitext, Springer, New York, 2012.
		
		
		\bibitem{CK}
		P. G. Casazza and J. Kova\v{c}evi\'{c}, Equal-norm tight frames with erasures, \emph{Adv. Comput. Math.}, 18 (2003), 387--430.
		
		
		\bibitem{CK1}
		P. G. Casazza and G. Kutyniok, Finite Frames: Theory and Applications, Birkh$\ddot{a}$user, 2012.
			
		\bibitem{C2}  
		O. Christensen, An introduction to frames and Riesz bases, Second edition,  Birkh$\ddot{a}$user,  2016.
		
		
		\bibitem{CS}
		D. Cvetkovi\'{c} and  S. Simi\'{c}, Graph spectra in computer science, \emph{Linear Algebra Appl.}, 434 (6) (2011), 1545--1562.
		
		
		\bibitem{DGM}
		I. Daubechies, A. Grossmann and Y. Meyer, Painless nonorthogonal expansions, \emph{J. Math. Phys.}, 27 (1986), 1271--1283.
		
		
		\bibitem{DA}
		Deepshikha and A. Samanta, Averaged numerically optimal dual frames for erasures,  \emph{Linear Multilinear Algebra},  71 (2) (2023),  301--316.
		
		\bibitem{DV}
		Deepshikha and L. K. Vashisht, A note on discrete frames of translates in $\mathbb{C}^N$,  \emph{TWMS  J.  Appl. Eng. Math.},  6 (1) (2016),  143--149.
		
		\bibitem{DV1}
		Deepshikha and L. K. Vashisht, Necessary and sufficient conditions for discrete wavelet frames in $\mathbb{C}^N$,  \emph{J. Geom. Phys.},  117 (2017),  134--143.
		
		\bibitem{DS}
		R. J. Duffin and A. C. Schaeffer, A class of nonharmonic Fourier series, \emph{Trans. Amer. Math. Soc.}, 72 (1952), 341--366.
		
		\bibitem{HP}
		R. Holmes and V. Paulsen, Optimal frames for erasures, \emph{Linear Algebra Appl.}, 377 (2004), 31--51.
		
		\bibitem{KOPT}
		G. Kutyniok, K. A. Okoudjou, F. Philipp and E. K. Tuley, Scalable frames, \emph{Linear Algebra Appl.}, 438 (5) (2013), 2225--2238.
		
		\bibitem{LH1}
		J. Lopez and D. Han, Optimal dual frames for erasures, \emph{Linear Algebra Appl.}, 432 (1) (2010), 471--482.
		
		\bibitem{MSS1}
		Adam W. Marcus, Daniel A. Spielman and N. Srivastava, Interlacing families I: Bipartite Ramanujan graphs of all degrees, \emph{Ann. Math.}, 182 (1) (2015), 307--325. 
		
		\bibitem{MSS}
		A. W. Marcus, D. A. Spielman, and N. Srivastava. Interlacing families II: Mixed characteristic polynomials and the Kadison—Singer problem, \emph{Ann. Math.},  182 (1) (2015), 327--350.
		
		
		\bibitem{PHM}
		S. Pehlivan, D. Han and R. Mohapatra, Linearly connected sequences and spectrally optimal dual frames for erasures, \emph{J. Funct. Anal.}, 265 (11) (2013), 2855--2876.
		
		
		\bibitem{SH}
		T. Strohmer and R.  W. Heath, Grassmannian frames with applications to coding and communication, \emph{Appl. Comput. Harmon. Anal.}, 14 (3) (2003), 257--275.
		
		
		
	\end{thebibliography}
\end{document}